\documentclass[fleqn,oneside,a4paper,12pt]{article}
\usepackage[english]{babel}
\usepackage[T1]{fontenc}
\usepackage[a4paper,tmargin=3cm,bmargin=3cm,rmargin=2.1cm,lmargin=2.1cm]{geometry}
\usepackage{amsfonts}
\usepackage{qsymbols}
\usepackage{amsmath}
\usepackage{euscript}
\usepackage{wasysym}
\usepackage{babel}
\usepackage{amssymb}
\usepackage[T1]{fontenc}
\usepackage[all]{xy}
\usepackage{graphicx}
\usepackage{epsf}
\usepackage{phonetic}
\usepackage{url}
\usepackage{showkeys}
\usepackage{stmaryrd}
\usepackage{mathrsfs}
\usepackage{expdlist}
\usepackage{fancyhdr}
\usepackage{theorem,amsfonts}
\theoremstyle{break}

\newtheorem{defi}{Definition}[section]
\newtheorem{exem}[defi]{Example}

\newtheorem{rem}[defi]{Remark}

\newtheorem{lem}[defi]{Lemma}
\newtheorem{theo}[defi]{Theorem}

\title{Friezes of type $\mathbb{D}$}
\author{Kodjo Essonana Magnani}
\begin{document}
\maketitle
\begin{abstract}
In this article, we establish a link between the values of a frieze of type $\mathbb{D}_{n}$ and some values of a particular frieze of type $\mathbb{A}_{2n-1}$. This link allows us to compute, independently of each other, all the cluster variables in the cluster algebra associated with a quiver $Q$ of type $\mathbb{D}_{n}.$ 
\end{abstract}

\section{Introduction}

\hspace{0.7 cm} Cluster algebras were introduced by S. Fomin and A. Zelevinsky in [FZ02, FZ03]. They are a class of commutative algebras which was shown to be connected to various areas of mathematics like combinatorics, Lie theory, Poisson geometry, Teichmüller theory, mathematical physics and representation theory of algebras.

\hspace{0.3 cm} A cluster algebra is generated by a set of variables, called \textit{cluster variables}, obtained recursively by a combinatorial process known as \textit{mutation} starting from a set of initial cluster variables. Explicit computation of cluster variables is difficult and has been extensively studied, see [ARS10, ADSS11, AD11, AR12, BMRRT06, BMR08].

\hspace{0.3 cm} In order to compute cluster variables, one may use \textit{friezes}, which were introduced by Coxeter [C71] and studied by Conway and Coxeter [CC73-I,II]. Various relationships are known between friezes and cluster algebras, see [CC06, M11, Pr08, ARS10, BM12, D10, AD11, ADSS12, AR12].

\hspace{0.3 cm} The present work is motivated by the use of friezes to compute cluster variables and is inspired by the result in [ARS10] giving an explicit formula as a product of $2 \times 2$ matrices for all cluster variables in coefficient-free cluster algebras of type $\mathbb{A}$, thus explaining at the same time the Laurent phenomenon and positivity.

\hspace{0.3 cm} Our objective here is to show that the same technique can be used for computing cluster variables in coefficient-free cluster algebras of type $\mathbb{D}$. The friezes of type $\mathbb{D}$ have already been studied by K. Baur and R. Marsh in [BM12] but our approach here is different.

\hspace{0.3 cm} In this paper, we establish a link between the values of a frieze of type $\mathbb{D}$ and some values of a particular frieze of type $\mathbb{A}$. For this we associate with each quiver of type $\mathbb{D}_{n}$ a particular quiver of type $\mathbb{A}_{2n-1}$. This correspondence allows us to use the algorithm of [ARS10] for computing cluster variables of a cluster algebra of type $\mathbb{D}_{n}$.

\hspace{0.3 cm} The article is organized as follows. In section 2, we recall some basic notions on friezes of type $\mathbb{A}$, set the preliminaries on friezes of type $\mathbb{D}$ and establish a correspondence between a frieze of type $\mathbb{D}$ and a particular frieze of type $\mathbb{A}$. In section 3, we give an algorithm to compute all the cluster variables in a cluster algebra associated with a quiver $Q$ of type $\mathbb{D}$. 

\section{Relation between friezes of types $\mathbb{A}$ and $\mathbb{D}$  }

In this section we establish a link between values of a frieze of type  $\mathbb{D}$ and some values of a particular frieze of type $\mathbb{A}$. For this we recall briefly some properties of the friezes of types $\mathbb{A}$ and $\mathbb{D}$.

\subsection{Some notions on friezes of type $\mathbb{A}$}

Let $n\geq 0$ be an integer. Consider a regular polygon $\mathcal{P}_{n+3}$ with $(n+3)$ sides such that its vertices are labelled $P_{1}, P_{2}, ..., P_{n+3}$ in clockwise orientation and let $\mathcal{T}$ be a triangulation of $\mathcal{P}_{n+3}$. An arc in $\mathcal{T}$ is the isotopy class of a curve joining two vertices $P_{i}, P_{j}$ of $\mathcal{P}_{n+3}$ with $j \notin \{i-1, \; i+1 \}$ and the convention that $P_{1} = P_{n+4}$. A triangulation is a maximal set of arcs which do not cross each other. With every triangulation  $\mathcal{T}$ of $\mathcal{P}_{n+3}$, we can associate a quiver $\Theta_{\mathcal{T}}$ in the following way.

\hspace{0.3 cm} The points of the quiver $\Theta_{\mathcal{T}}$ represent arcs of the triangulation $\mathcal{T}$. There is an arrow $\alpha \rightarrow \beta$ if there exists a triangle $\Delta$ of $\mathcal{T}$ such that $\alpha$ and $\beta$ are two sides of $\Delta$ and if $\beta$ is successor of $\alpha$ in $\Delta$ with respect to clockwise orientation about the common vertex of $\alpha$ and $\beta$. Thus the quiver $\Theta_{\mathcal{T}}$ has $n$ points, $n$ being the number of arcs in the triangulation $\mathcal{T}$ associated with $\mathcal{P}_{n+3}$ (see [CCS06-3.2] of Caldero, Chapoton and Schiffler and [FST08-4.1] of Fomin, Shapiro and Thurston).

\hspace{0.3 cm} In this paper we consider only triangulations without internal triangles (a triangle is call internal if none of its sides coincides with a side of the polygon). A quiver associated with a triangulation of $\mathcal{P}_{n+3}$  has therefore underlying graph of the Dynkin type $\mathbb{A}_{n}$.

\hspace{0.3 cm} Let $\Theta$ be a finite acyclic (containing no oriented cycles) quiver with $\Theta_{0}$ the set of its points and $\Theta_{1}$ the set of its arrows and $\mathbb{K}$ a field. The translation quiver $\mathbb{Z}\Theta$ associated with $\Theta$  (see [ASS-VIII.1.1]) is consisting of two sets: the set of points  $(\mathbb{Z}\Theta)_{0} = \mathbb{Z} \times \Theta_{0} = \left\lbrace  (k, i) | k \in \mathbb{Z}, \, i \in \Theta_{0} \right\rbrace $ and the set of arrows\\
 $(\mathbb{Z}\Theta)_{1} = \left\lbrace  (k, \alpha): (k, i)\rightarrow (k, j) | k \in \mathbb{Z}, \, \alpha: i \rightarrow j \in \Theta_{1} \right\rbrace  \cup \\ \left\lbrace  (k, \alpha'): (k, j)\rightarrow (k+1, i) | k \in \mathbb{Z}, \, \alpha: i \rightarrow j \in \Theta_{1} \right\rbrace. $ Let us define a frieze associated with the quiver $\Theta$.

\hspace{0.3 cm} In the translation quiver $\mathbb{Z}\Theta$, let us replace the points $(k, i) \in (\mathbb{Z}\Theta)_{0}$ by their images $\mathfrak{a}(k, i)$ obtained applying a frieze function $\mathfrak{a}$: $(\mathbb{Z}\Theta)_{0} \rightarrow \mathbb{K}$ defined for some initial values $\mathfrak{a}(0,i) \in \mathbb{K}$ as follows:\\ $\mathfrak{a}(k,i)\mathfrak{a}(k+1,i) = 1+ \prod_{(k,i)\rightarrow (m,j)} \mathfrak{a}(m,j)$ where the product is taken over the arrows (see [ARS10-2]). The resulting translation quiver with values associated with its vertices is called a frieze.\\

\hspace{0.3 cm} We now recall the definition of seed due to Fomin and Zelevinsky [FZ03-1.2].\\

\begin{defi}
 Let $\Gamma$ be a quiver with $n$ points and $\chi = \{u_{1},u_{2}, ..., u_{n}\}$ a set of variables called cluster variables, such that a variable $u_{i}$ is associated with the point $i$ (with $1\leq i \leq n$) of $\Gamma$. The set $\chi$ is called a cluster and the pair $(\Gamma , \chi)$ is called a seed.\\
\end{defi}

 \hspace{0.3 cm} We can obtain other seeds by mutation (see [FZ03-1.2]) starting from the seed $(\Gamma, \chi)$. The set of all cluster variables obtained by successive mutation generates an algebra over $\mathbb{Z}$ called \textit{cluster algebra} which is denoted by $\mathcal{A}(\Gamma, \chi)$.
 
\begin{rem}
We say that the type of a seed $(\Gamma, \chi)$ and cluster algebra $\mathcal{A}(\Gamma, \chi)$ coincides with the type of quiver $\Gamma$.
\end{rem}

\hspace{0.3 cm} Note that if we associate with each vertex $i$ of $\Theta$ a variable $u_{i}$ and take initial values $\mathfrak{a}(o, i ) = u_{i}$ then all the values of the frieze associated to $\Theta$ are cluster variables of the cluster algebra $\mathcal{A}(\Theta, \chi)$ (see [AD11]).\\

\hspace{0.3 cm} Our aim in this paper is to compute cluster variables of a cluster algebra of type $\mathbb{D}_{n}$. We start by recalling some properties of friezes of type $A_{n}$ which will be used in our construction.

\begin{rem}
In what follows, we assume that in the case of quiver $\Theta$ of type $\mathbb{A}_{n}$ all arrows of $\mathbb{Z}\Theta$ are directed either south-east or north-east (if all arrows of $\Theta$ are drawn horizontally, the arrows oriented from left to right in $\Theta$ become directed north-east in $\mathbb{Z}\Theta$ and the arrows oriented from right to left in $\Theta$ are directed south-east in $\mathbb{Z}\Theta$). When talking about friezes, we identify a point $(k,i) \in (\mathbb{Z}\Theta)_{0}$ with its image $\mathfrak{a}(k,i)$.\\ 

\end{rem}

\hspace{0.3 cm} Now we define the notion of a diagonal in a frieze of type $\mathbb{A}_{n}$.

\begin{defi}
Let $k_{0} \in \mathbb{Z}$ and $\Theta$ be a finite acyclic quiver of type $\mathbb{A}_{n}$. A descending (or ascending) diagonal in the frieze associated with $\Theta$ is the concatenation of $(n-1)$ arrows with south-east (or north-east) orientation starting at $\mathfrak{a}(k_{0},n)$ (or $\mathfrak{a}(k_{0}, 1)$, respectively).\\

\end{defi}

\hspace{0.3 cm} Let us introduce the notion of a fundamental quiver in a frieze of type $\mathbb{A}_{n}$. \\ Fix an integer $k_{0}$. Let $d_{1}^{k_{0}}$ be the descending diagonal starting at $(k_{0},n) \in (\mathbb{Z}\Theta)_{0}$ and $d_{2}^{k_{0}}$ the ascending diagonal ending at $(k_{0}+n,n) \in (\mathbb{Z}\Theta)_{0}$.

\begin{defi}
Let $\Theta$ be a finite acyclic quiver of  type $\mathbb{A}_{n}$. The fundamental quiver $\Theta_{f}^{k_{0}}$ in the frieze associated with $\Theta$ is the portion of the translation quiver $\mathbb{Z}\Theta$ bounded by diagonals $d_{1}^{k_{0}}$ and $d_{2}^{k_{0}}$ inclusively.\\
\end{defi}

\hspace{0.3 cm} Independently of the chosen integer $k_{0}$, the values of the frieze function associated with vertices of the fundamental quiver $\Theta_{f}^{k_{0}}$ form the set of all cluster variables of the associated cluster algebra $\mathcal{A}(\Theta, \chi)$ ([FZ03]). In what follows, we omit the index $k_{0}$ and specify, where necessary, the position of diagonals $d_{1}$, $d_{2}$ in the translation quiver in an alternative way.

\hspace{0.3 cm} We give an example of a fundamental quiver in the frieze associated with  $\Theta$ of type $\mathbb{A}_{3}$.

\begin{exem}
Consider the following quiver $\Theta$ of type $\mathbb{A}_{3}$ : $ 1\rightarrow 2\rightarrow 3$.\\
The translation quiver $\mathbb{Z}\Theta$ associated with $\Theta$ is:\\
$
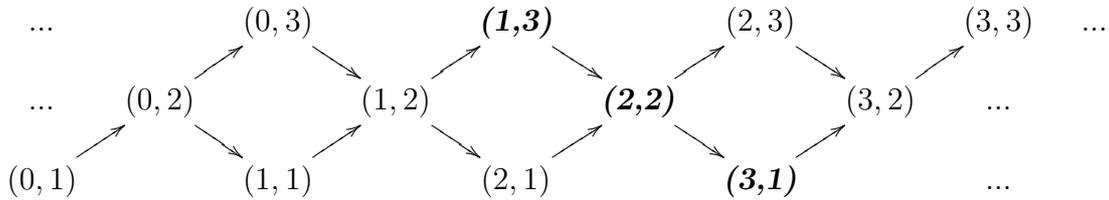
\xymatrix@=1pc{
&   & &  & &   & &   & &   & &  &&      &  \\
 & ...  &   &(0,3)\ar[dr]  & & \textbf{(1,3)}\ar[dr] & & (2,3)\ar[dr] &&  (3,3) &... &  &  \\
 &...& (0,2)\ar[dr]\ar[ur] & & (1,2)\ar[ur]\ar[dr] & & \textbf{(2,2)}\ar[dr]\ar[ur]  & &(3,2)\ar[ur] &... & &  & &        \\
   &  (0,1)\ar[ur] &   & (1,1)\ar[ur] && (2,1)\ar[ur] & &  \textbf{(3,1)}\ar[ur] &&...& & & &  &   \\
&    &  &  & &    & &   & &    &&   & &   &           \\ 
}$\\

Associating a value $\mathfrak{a}(k, i)$ of the frieze function with every vertex we obtain an example of a frieze.\\

A fundamental quiver is (with $k_{0} = 0$):\\

$
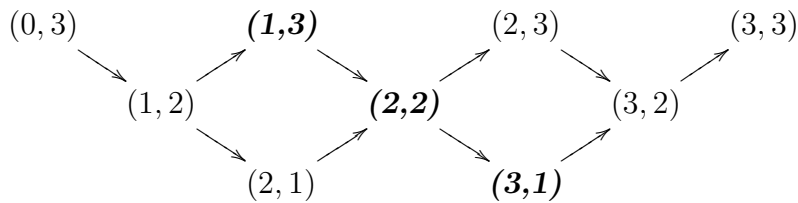
\xymatrix@=1pc{
&   & &  & &   & &   & &   & &  &&      &  \\
 &  &   &(0,3)\ar[dr]  & &\textbf{(1,3)}\ar[dr] & &(2,3)\ar[dr] &&(3,3) & &  &  \\
 &&  & &(1,2)\ar[ur]\ar[dr] & &\textbf{(2,2)}\ar[dr]\ar[ur]  & &(3,2)\ar[ur] & & &  & &        \\
   &   &   &  &&(2,1)\ar[ur] & & \textbf{(3,1)}\ar[ur] &&& & & &  &   \\
&    &  &  & &    & &   & &    &&   & &   &           \\ 
}$\\
\end{exem}

\begin{rem}
As a consequence of the definition of the frieze function, all the squares of the form

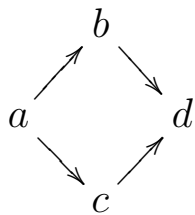
\resizebox{7cm}{!}{
$\xymatrix@=1pc{ 
 &   & &  b\ar[dr]  & &&&        \\
   &\text{}   & a\ar[ur]\ar[dr]   &      & d     & &&        \\
 &   &    &   c\ar[ur]&    &      &      &    \\
}$}\\  in the frieze of type $\mathbb{A}$ satisfy the relation $ad-bc = 1$, with $a, b, c, d \in \mathbb{K}$. This relation is called uni-modular rule.
\end{rem}

\subsection{Preliminaries on friezes of type $\mathbb{D}$}

Consider a quiver $Q$ whose underlying graph is of type $\mathbb{D}_{n}$. We agree to label the points of $Q$ as follows:  $
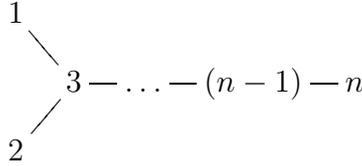
\xymatrix @=10pt
{
&1
&
&
&
&&\\
\text{}
&
&3 \ar@{-}[r]\ar@{-}[ul]\ar@{-}[dl]
&\ldots \ar@{-}[r]
&(n-1) \ar@{-}[r]
&n \text{} \\
&2
&
&
}
$\\ 
where a solid segment represents an arrow without its orientation.\\
The \textit{fork} is the full sub-quiver of $Q$ generated by the points $\{1, 2, 3\}$.\\
We agree to call:\\
- \textit{fork arrows}, the arrows of the fork, \\
- \textit{joint} of the fork, the point $3$,\\
- and \textit{fork vertices}, the points $1$ and $2$.\\

\hspace{0.3 cm} If we associate with each vertex $i$ of $Q$ a variable $u_{i}$ then we get the seed $\mathcal{G}= (Q , \chi)$ whose underlying graph can be represented by the following diagram: $$
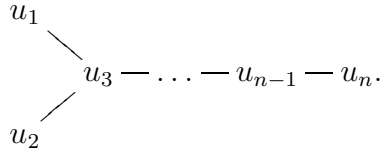
\xymatrix @=10pt
{
&u_{1}
&
&
&
&&\\
\text{}
&
&u_{3} \ar@{-}[r]\ar@{-}[ul]\ar@{-}[dl]
&\ldots \ar@{-}[r]
&u_{n-1}\ar@{-}[r]
&u_{n} \text{} .\\
&u_{2}
&
&
}
$$\\ 

\hspace{0.3 cm} Let us define a frieze on $\mathbb{Z}\mathbb{D}_{n}$ by the function $\mathfrak{a}$: $(\mathbb{Z}\mathbb{D}_{n})_{0} \rightarrow \mathbb{Q}(u_{1}, u_{2},...,u_{n})$ such that for $(k,i) \in (\mathbb{Z}\mathbb{D}_{n})_{0}$ we have: $\mathfrak{a}(k,i)\mathfrak{a}(k+1,i) = 1+ \prod_{(k,i)\rightarrow (m,j)} \mathfrak{a}(m,j)$ with initial variables $\mathfrak{a}(0, i) = u_{i}$. Then, all the values of the frieze are cluster variables and all cluster variables are represented in the frieze (see [AD11] of Assem and Dupont).\\
 
\hspace{0.3 cm} Consider a seed $\mathcal{G}$ of type $\mathbb{D}_{n}$. Let $\mathcal{F}$ be the part of the frieze associated with $\mathcal{G}$, which is given by the set of vertices $\{ (k , i) \in (\mathbb{Z}\mathbb{D}_{n})_{0} \,| k = 0,1,...,n; i = 1,...,n \}$ and their incident arrows. $\mathcal{F}$ contains all cluster variables of cluster algebra $\mathcal{A}(\mathcal{G})$ of type $\mathbb{D}_{n}$ ([FZ03]) and a second copy of initial variables $\{u_{i}\}_{i=1}^{n}$. \\

\hspace{0.3 cm} We give in the following example the part $\mathcal{F}$ corresponding to various seeds of type $\mathbb{D}_{5}$.\\

\begin{exem}
\begin{enumerate}

\item $\xymatrix @=10pt{
&u_{1}\ar[dr]
&
&
&\text{}\\
\text{For the following seed of type $\mathbb{D}_{5}$ :}\quad
&
&u_{3}\ar[r]
&u_{4}\ar[r]
&u_{5}
&&\text{}\\  
&u_{2}\ar[ur]
&
&
}  $\\

 $\mathcal{F}$ has the form:\\
\resizebox{12cm}{!}{
$\xymatrix@=1pc{
& &   &   &    &   &    &      &      &&&&&&&&         \\
& &&   &  &  u_{5}\ar[dr] && .\ar[dr]&  & .\ar[dr] &   & .\ar[dr] &  &.\ar[dr]&&u_{5}& && \\
&&&   &u_{4}\ar[dr]\ar[ur] & & .\ar[dr]\ar[ur]& &  .\ar[dr]\ar[ur]&   &.\ar[ur]\ar[dr]&&.\ar[dr]\ar[ur]&      &u_{4}\ar[ur]&&\\
&&&u_{3} \ar[dr]\ar[ddr]\ar[ur]&&.\ar[ddr]\ar[ur]\ar[dr]& &.\ar[dr]\ar[ddr]\ar[ur]& &   .\ar[dr]\ar[ddr]\ar[ur]& & .\ar[dr]\ar[ddr]\ar[ur] &          &u_{3}\ar[ur]&&&&&\\
&& u_{1} \ar[ur]&&\frac{(1+u_{3})}{u_{1}}\ar[ur] &  &.\ar[ur]&  & .\ar[ur] & &  .\ar[ur] &   &u_{2}\ar[ur]&&&&      \\
&&u_{2} \ar[uur]&&\frac{(1+u_{3})}{u_{2}}\ar[uur]&&.\ar[uur] & &.\ar[uur]&&.\ar[uur]& &u_{1}\ar[uur]&&&&&&\\
 & &   &  &    &   &      &    &&&&&&&\\
}$} \\

\item  $\xymatrix @=10pt{
&u_{1}
&
&
&\text{}\\
\text{For the following seed of type $\mathbb{D}_{5}$ :}\quad
&
&u_{3}\ar[r]\ar[ul]
&u_{4}\ar[r]
&u_{5}
&&\text{}\\  
&u_{2}\ar[ur]
&
&
}  $\\

$\mathcal{F}$ has the form:\\
\resizebox{12cm}{!}{
$\xymatrix@=1pc{
& &   &   &    &   &    &      &      &&&&&&&&         \\
& &&   &  &  u_{5}\ar[dr] && .\ar[dr]&  & .\ar[dr] &   & .\ar[dr] &  &.\ar[dr]&&u_{5}& && \\
&&&   &u_{4}\ar[dr]\ar[ur] & & .\ar[dr]\ar[ur]& &  .\ar[dr]\ar[ur]&   &.\ar[ur]\ar[dr]&&.\ar[dr]\ar[ur]&      &u_{4}\ar[ur]&&\\
&&&u_{3} \ar[dr]\ar[r]\ar[ur]&u_{1}\ar[r]&.\ar[r]\ar[ur]\ar[dr]&.\ar[r] &.\ar[dr]\ar[r]\ar[ur]&.\ar[r] &   .\ar[dr]\ar[r]\ar[ur]&.\ar[r] & .\ar[dr]\ar[r]\ar[ur] &u_{2}\ar[r]      &u_{3}\ar[ur]\ar[rd]&&&&&\\
&& u_{2} \ar[ur]&&\frac{(1+u_{3})}{u_{2}}\ar[ur] &  &.\ar[ur]&  & .\ar[ur] & &  .\ar[ur] &   &\frac{(1+u_{3})}{u_{1}}\ar[ur]&&u_{1}&&      \\
&&&&&& & &&&&&&&&&&\\
 & &   &  &    &   &      &    &&&&&&&\\
}$} \\
\end{enumerate}
\end{exem}

The study of a seed of type $\mathbb{D}_{n}$ can be subdivided into three cases depending on the orientation of fork arrows:\\
- The fork is composed by two arrows leaving the joint,\\
- The fork is composed by two arrows entering the joint,\\
- The fork is composed by one arrow leaving the joint and another arrow entering the joint.\\

The following lemma which is a reformulation of lemma $5.2$ in [ASS-VII] of Assem, Simson and Skowro\'nski allows us to reduce the study to one of the three cases. We refer to [FZ03-1.2] and [FZ03-8] for notions of mutation and mutation equivalent, respectively.

\begin{lem}
Let $Q_{1}$ and $Q_{2}$ be two quivers having the same underlying graph $G$. If $G$ is a tree then $Q_{1}$ and $Q_{2}$ are mutation equivalent. $\square$ \\
\end{lem}

\hspace{0.3 cm} Consider two seeds $\mathcal{G}_{1}$ and  $\mathcal{G}_{2}$ of type $\mathbb{D}_{n}$ which are transformed into one another by a mutation on a fork vertex corresponding to the variable $\upsilon$ of $\mathcal{G}_{1}$ and $\tilde{\upsilon}$ of $\mathcal{G}_{2}$.

\hspace{0.3 cm} Suppose that the quiver $Q_{1}$ associated with $\mathcal{G}_{1}$ has a fork consisting of two arrows leaving (or entering) the joint and the quiver $Q_{2}$ associated with $\mathcal{G}_{2}$ has a fork consisting of one arrow leaving and another entering.

\hspace{0.3 cm} The seeds $\mathcal{G}_{1}$ and  $\mathcal{G}_{2}$ generate the same cluster algebra $\mathcal{A}.$ Therefore each seed generates by mutations the set of all cluster variables of $\mathcal{A}.$ It is known that $\mathbb{Z}Q$ is periodic if $Q$ is of Dynkin type and the corresponding part $\mathcal{F}$ contains all the cluster variables of $\mathcal{A}$. According to Proposition $2.1$ in [S08] of Schiffler, also in [H88-I.5.6] of Happel, we deduce that the parts $\mathcal{F}_{1}$ corresponding to $\mathcal{G}_{1}$ and  $\mathcal{F}_{2}$ corresponding to $\mathcal{G}_{2}$ contain the same and all cluster variables of the cluster algebra  $\mathcal{A}$.

\hspace{0.3 cm} Since our aim is to compute these cluster variables independently of each other, it suffices to study only the case where the fork consists of two arrows entering  (or leaving) the joint.\\

In the following, $\mathcal{G}$ denotes a seed of type $\mathbb{D}_{n}$ with a quiver whose fork is composed by two arrows entering (or leaving) the joint.\\

\begin{defi}
We call modelled quiver $\bar{\mathcal{F}}$ associated with $\mathcal{G}$, the part of a translation quiver obtained from $\mathcal{F}$ as follows:\\
$1.$ by gluing in $\mathcal{F}$ the arrows of each shifted copy of the fork,\\
$2.$ by multiplying the values of vertices of the fork corresponding to the arrows that were glued in step $1$. Namely the arrows obtained by gluing in step $1$ have $\mathfrak{a}(k, 1)\mathfrak{a}(k, 2)$ as the corresponding variables.\\
\end{defi}

We give an example of a modelled quiver associated with a seed of type $\mathbb{D}_{5}$.

\begin{exem}

Consider the item $1.$ of example $2.8$.\\
The modelled quiver $\bar{\mathcal{F}}$ obtained from $\mathcal{F}$ has the form:\\
\resizebox{12cm}{!}{
$
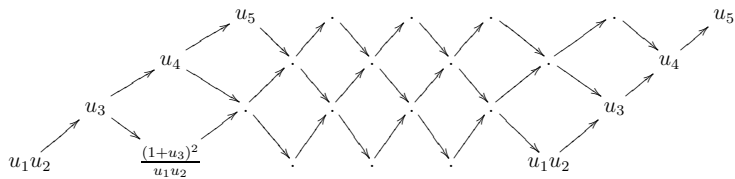
\xymatrix@=1pc{
& &   &   &    &   &    &      &      &&&&&&&&         \\
& &&   &  &  u_{5}\ar[dr] && .\ar[dr]&  & .\ar[dr] &   & .\ar[dr] &  &.\ar[dr]&&u_{5}& && \\
&&&   &u_{4}\ar[dr]\ar[ur] & & .\ar[dr]\ar[ur]& &  .\ar[dr]\ar[ur]&   &.\ar[ur]\ar[dr]&&.\ar[dr]\ar[ur]&      &u_{4}\ar[ur]&&\\
&&&u_{3} \ar[dr]\ar[ur]&&.\ar[ur]\ar[dr]& &.\ar[dr]\ar[ur]& &   .\ar[dr]\ar[ur]& & .\ar[dr]\ar[ur] &          &u_{3}\ar[ur]&&&&&\\
&& u_{1}u_{2}\ar[ur]&&\frac{(1+u_{3})^{2}}{ u_{1}u_{2}}\ar[ur] &  &.\ar[ur]&  & .\ar[ur] & &  .\ar[ur] &   &u_{1}u_{2}\ar[ur]&&&&      \\
 & &   &  &    &   &      &    &&&&&&&\\
}$}\\

\end{exem}

\subsection{Correspondence between a frieze of type $\mathbb{D}$ and a particular frieze of type $\mathbb{A}$}

We let $\bar{\omega}$ be the full sub-quiver of $Q$ generated by all points except the point $2$, and drawn from left to right in such a way that the vertices appear in increasing order. Note that \textbf{$\bar{\omega}$} is a quiver of type $\mathbb{A}_{n-1}$.  This leads to a way to associate a quiver of type $\mathbb{D}$ with a quiver of type $\mathbb{A}$, which, in turn, will allow us to associate a frieze of type $\mathbb{D}$ with a particular frieze of type $\mathbb{A}$. 

\hspace{0.3 cm} For a quiver $\Theta$ of type $\mathbb{A}_{n}$ (drawn from left to right), let us denote by $^{t}\Theta$ the transpose of $\Theta$, that is, the quiver obtained by redrawing $\Theta$ from right to left.\\

 With a quiver $Q$ of type $\mathbb{D}_{n}$ we associate the quiver $Q'$ of the form\\ $^{t}\bar{\omega}\rightarrow 0 \rightarrow \bar{\omega}$. Note that the underlying graph of $Q'$ is of type $\mathbb{A}_{2n-1}$ labelled as follows:\\
 $\xymatrix @=10pt
{
&
&
&&
&& \\
\text{}
&n \ar@{-}[r]
&(n-1) \ar@{-}[r] 
&\ldots \ar@{-}[r] 4
&3 \ar@{-}[r] 
&1 \ar@{-}[r]
&0 \ar@{-}[r]
&1 \ar@{-}[r]
&3 \ar@{-}[r]
&4 \ldots
&(n-1) \ar@{-}[l]
&n. \ar@{-}[l]
& 
}
$\\

The following example shows how to construct the quiver $Q'$ from a given quiver $Q$ of type $\mathbb{D}_{4}$.\\

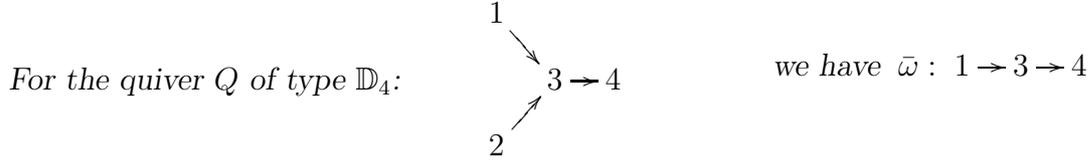
\begin{exem}
 $\xymatrix @=10pt{
&{1}\ar[dr]
&
&
&\text{}\\
\text{For the quiver $Q$ of type $\mathbb{D}_{4}$: }\quad
&
&{3}\ar[r]
&{4}
&&\text{}\\  
&{2}\ar[ur]
&
&
}  $
$\xymatrix @=10pt{
&
&
&\text{}\\
& \text{we have} \; \; \textbf{$\bar{\omega}$} :\;{1}\ar[r]
&{3}\ar[r]
&4
&&\text{}\\  
&
&
&
&
}  $

$\xymatrix @=10pt{
&
&
&\text{}\\
& \text{and} \; \; \textbf{${Q'}$} :\; {4}
&{3}\ar[l]
&{1}\ar[l]\ar[r]
&{0}\ar[r]
&{1}\ar[r]
&{3}\ar[r]
&{4}.
&&\text{}\\  
&
&
&
&
}  $\\

\end{exem}

Now we associate variables with the points of $Q'$ as follows:\\
$1.$ with the vertex labelled by $i=1$ we associate the product $u_{1}u_{2}$,\\
$2.$ with a vertex labelled by $i \neq 1$ we associate the variable $u_{i}$. (Note that by doing so we introduce a new variable $u_{0}$ which is not a cluster variable of $\mathcal{G}$. This variable will be later replaced by $1$).\\

The result of this association is a seed which will be denoted by $\Lambda'$. Its underlying graph with variables is:\\$$\xymatrix @=10pt
{
&
&
&&
&& \\
\text{}
&u_{n} \ar@{-}[r]
&u_{n-1} \ar@{-}[r] 
&\ldots \ar@{-}[r] u_{4}
&u_{3} \ar@{-}[r] 
&u_{1}u_{2} \ar@{-}[r]
&u_{0} \ar@{-}[r]
&u_{1}u_{2} \ar@{-}[r]
&u_{3} \ar@{-}[r]
&u_{4} \ldots
&u_{n-1} \ar@{-}[l]
&u_{n}. \quad (1) \ar@{-}[l]
& 
}
$$\\

\hspace{0.3 cm} The seed $\Lambda'$ can be represented as follows: $^{t}\bar{\omega}_{*}\rightarrow u_{0} \rightarrow \bar{\omega}_{*}$, where $\bar{\omega}_{*}$ is $\bar{\omega}$ with variables.

\hspace{0.3 cm} The graph underlying the seed $\Lambda'$ contains $(2n-1)$ vertices. We enumerate them from left to right in $(1)$ so that the variable $u_{1}u_{2}$ is associated to $(n-1)$th and $(n+1)$th vertices.

\hspace{0.3 cm} Since $u_{0}$ is not a cluster variable of $\mathcal{G}$, in the rest of this paper we evaluate  $u_{0}= 1$.

\hspace{0.3 cm} The following theorem establishes a link between the modelled quiver $\bar{\mathcal{F}}$ associated with a seed $\mathcal{G}$ of type $\mathbb{D}_{n}$ and a fundamental quiver in the frieze associated with $\Lambda'$.\\

\begin{theo} 
Let $n\geq 4$ be an integer, $\mathcal{G}$ a seed of type $\mathbb{D}_{n}$ and $\bar{\mathcal{F}}$ the modelled quiver associated with $\mathcal{G}$.\\
Then $\bar{\mathcal{F}}$ is a full sub-quiver of a fundamental quiver in the frieze of type $\mathbb{A}_{2n-1}$ associated with $\Lambda'$. The descending diagonal $d_{1}$ which bounds this fundamental quiver passes through the $(n+1)$th vertex of $\Lambda'$ ( whose associated variable is $u_{1}u_{2}$ ).\\
\end{theo}

\textbf{Proof:} We give the proof for a seed $\mathcal{G}$ whose associated quiver has two arrows entering the fork. The case of a quiver whose fork arrows leave the joint is absolutely analogous.  Let $\bar{\mathcal{F}}$ be the modelled quiver associated with ${\mathcal{G}}$. All the squares of the form

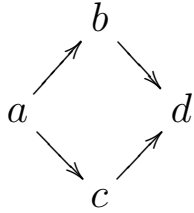
\resizebox{7cm}{!}{
$\xymatrix@=1pc{ 
 &   & &  b\ar[dr]  & &&&        \\
   &\text{}   & a\ar[ur]\ar[dr]   &      & d     & &&        \\
 &   &    &   c\ar[ur]&    &      &      &    \\
}$}\\ in $\bar{\mathcal{F}}$ satisfy the relation $ad-bc = 1$ called uni-modular rule, \\with $a, b, c, d \in \mathbb{Q}(u_{1}, u_{2},..., u_{n})$, as a consequence of the definition of a frieze.

\hspace{0.3 cm} We shall construct a frieze $\Upsilon$ of type $\mathbb{A}$ containing $\bar{\mathcal{F}}$. To this end, we shall complete the modelled quiver $\bar{\mathcal{F}}$ downward using the uni-modular rule:

\hspace{0.3 cm} The horizontal bottom line of $\bar{\mathcal{F}}$ which contains $u_{1}u_{2}$ is composed of $(n+1)$ values. By the uni-modular rule we can extend $\bar{\mathcal{F}}$ downward, thus creating the next line composed of $n$ values, the one after composed of  $(n-1)$ values and so on until the last line composed of one single entry. This entry will serve as the end point of the descending diagonal $d_{1}$.

\hspace{0.3 cm} Thus we have constructed $n$ new lines and we have $(2n-1)$ lines in total.\\
If we consider these $(2n-1)$ lines as part of a frieze $\Upsilon$ of type $\mathbb{A}_{2n-1}$, we can construct the entire frieze $\Upsilon$ by extending the existing pattern horizontally using the uni-modular rule. In particular we construct completely the descending diagonal $d_{1}$ and the ascending diagonal $d_{2}$ bounding the fundamental quiver of the frieze $\Upsilon$ of type $\mathbb{A}_{2n-1}$ .

\hspace{0.3 cm} Then we have constructed the fundamental quiver of the frieze $\Upsilon$ of type $\mathbb{A}_{2n-1}$ which contains $\bar{\mathcal{F}}$ as a full sub-quiver.\\

The next step is to prove that $\Lambda'$ generates the constructed frieze $\Upsilon$.\\

Taking into account that $\bar{\mathcal{F}}$ is contained in the fundamental quiver of the frieze $\Upsilon$ of type $\mathbb{A}_{2n-1}$, we deduce that there exists a quiver $\tilde{Q}$ of type $\mathbb{A}_{2n-1}$ which generates the frieze $\Upsilon$ such that $\tilde{Q}$ contains $\bar{\omega}$ (recall that $\bar{\omega}_{*}$ bounds $\bar{\mathcal{F}}$ on the right and on the left). Therefore $\tilde{Q}$ is of the form $\tilde{Q} = \tilde{Q}'-\bar{\omega}$  where the solid segment represents an arrow without its orientation and $\tilde{Q}'$ is an acyclic quiver of type $\mathbb{A}_{n}$. We have  $^{t}\tilde{Q} =\, ^{t}\bar{\omega}-\, ^{t}\tilde{Q}'$. The property of periodicity of the frieze of type $\mathbb{A}_{n}$ (see [CC73 I-II.21]) implies that the quivers $\tilde{Q}$ and $^{t}\tilde{Q}$ generate the same frieze. The periodicity of the frieze $\Upsilon$ and the presence of the second copy of $\bar{\omega}_{*}$ in $\bar{\mathcal{F}}$ implies that the quiver $^{t}\tilde{Q}$ contains $\bar{\omega}$ and is of the form $^{t}\tilde{Q} = \tilde{Q}''-\bar{\omega}$ where $\tilde{Q}''$ is an acyclic quiver of type $\mathbb{A}_{n}$. 

\hspace{0.3 cm} Because $^{t}\tilde{Q}$ can be written either as $^{t}\bar{\omega}-\, ^{t}\tilde{Q}'$ or as $\tilde{Q}''-\bar{\omega}$, and counting the number of points in $\tilde{Q}$, we deduce that $^{t}\tilde{Q}=\;  ^{t}\bar{\omega}-q-\bar{\omega}$   where $q$ is a point. Therefore $\tilde{Q}$ is of the form $ ^{t}\bar{\omega}-q-\bar{\omega}$ and the quiver $\tilde{Q}$ with associated variables is of the form $ ^{t}\bar{\omega}_{*}-u_{q}-\bar{\omega}_{*}$. \\

We are going to determine entirely the quiver $\tilde{Q}$ of type $\mathbb{A}_{2n-1}$, that is, to find the point $q$ and its incident arrows. To this end, we consider the following part of the frieze $\Upsilon$:\\

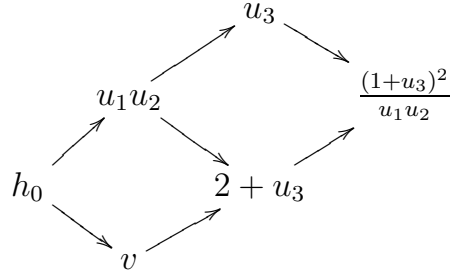
 \resizebox{13cm}{!}{
 $\xymatrix@=1pc{ 
&&& &   & &  u_{3}\ar[dr]  & &&&        \\
 &&&  &   & u_{1}u_{2}\ar[ur] \ar[dr]  &      & \frac{(1+u_{3})^{2}}{u_{1}u_{2}}    &&&        \\
&&&   &h_{0} \ar[ur]\ar[dr] & &   2+u_{3}\ar[ur]&    &      &   &   &    \\
 &&&&&v\ar[ur]&&&&&&&&\\
 &&&&&&&&&&&\\
}$}\\
where $h_{0}$ and $v$ are unknown and where $u_{1}u_{2}$ and $u_{3}$ correspond to points of $\bar{\omega}_{*}$.

\hspace{0.3 cm} Due to the form of the quiver $\tilde{Q}$ of type $\mathbb{A}_{2n-1}$, the point $q$ can be one of the two points underlying the values $h_{0}$ or $(2+u_{3})$. In other words, $u_{q}$ associated with the point $q$ may take two possible values $u_{q} = h_{0}$ or $u_{q} = (2+u_{3})$.

\hspace{0.3 cm} If we take all initial values $u_{k}=1$, the underlying graph of the quiver $\tilde{Q}$ with variables is of the form $\xymatrix @=10pt
{
\text{}
&1 \ar@{-}[r]
&1 \ar@{-}[r] 
&\ldots \ar@{-}[r] 
&1 \ar@{-}[r] 
&u_{q} \ar@{-}[r]
&1 \ar@{-}[r]
& \ldots
&1 \ar@{-}[l]
&1. \ar@{-}[l]
& 
}
$ In other words the value $u_{q}$ has two neighbours equal to one. According to the uni-modular rule, we have $\displaystyle v=\frac{h_{0}(2+u_{3})-1}{u_{1}u_{2}}$. If $u_{k} = 1$ for all $k$, the value $h_{0} \in \mathbb{Z}$, therefore $v \neq 1$. We conclude that the underlying vertex of the value $v$ is not in $\tilde{Q}$. Therefore the quiver $\tilde{Q}$ with associated variables is $ ^{t}\bar{\omega}_{*}\rightarrow h_{0} \rightarrow \bar{\omega}_{*}$ or $ ^{t}\bar{\omega}_{*}\leftarrow (2 + u_{3}) \leftarrow\bar{\omega}_{*}$ . \\

\hspace{0.3 cm} Suppose that $u_{q} = 2+u_{3}$.\\
 
 In this case, the form of $\tilde{Q}$ with variables ($ ^{t}\bar{\omega}_{*}\leftarrow u_{q} \leftarrow\bar{\omega}_{*}$) and the fact that the descending diagonal $d_{1}$ passes through the vertex of the variable $u_{1}u_{2}$ imply that $d_{1}$ contains both variables $u_{1}u_{2}$ and $u_{q} = 2+u_{3}$. This diagonal contains also the other variable $u_{1}u_{2}$ (according to the form of $\tilde{Q}$, the variable $u_{q}$ has two neighbours whose values are $u_{1}u_{2}$). According to the result in [C71-8] of Coxeter for a diagonal in a frieze, $(2+u_{3})$ divides the sum of values of its neighbours $(u_{1}u_{2}+u_{1}u_{2}) = 2u_{1}u_{2}$, that is under our assumption, $(2+u_{3})$ should divide $2u_{1}u_{2}$. This is false because for $u_{k} = 1$, \; the integer $3$ does not divide $2$. 

\hspace{0.3 cm} We conclude that $u_{q} = h_{0}$. Thus the quiver $\tilde{Q}$ is: $ ^{t}\bar{\omega}\rightarrow q \rightarrow \bar{\omega}$. 

\hspace{0.3 cm} Due to the fact that an arrow in $\bar{\omega}$ oriented from left to right becomes in $^{t}\bar{\omega}$ an arrow oriented from right to left, we have that $u_{1}u_{2}\rightarrow u_{3}$ in $\bar{\omega}_{*}$ becomes $u_{3}\leftarrow u_{1}u_{2}$ in $^{t}\bar{\omega}_{*}$. Then the frieze $\Upsilon$ contains the following diagram:\\
 \resizebox{13cm}{!}{
 $
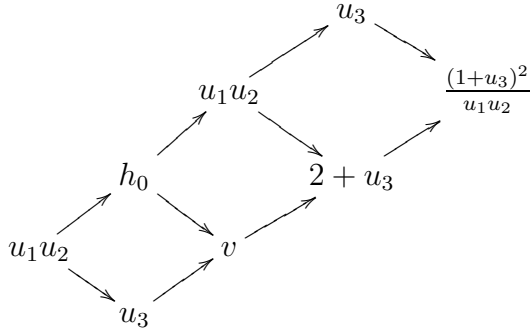
\xymatrix@=1pc{ 
&&& &   & &  u_{3}\ar[dr]  & &&&        \\
 &&&  &   & u_{1}u_{2}\ar[ur] \ar[dr]  &      & \frac{(1+u_{3})^{2}}{u_{1}u_{2}}    &&&        \\
&&&   &h_{0} \ar[ur]\ar[dr] & &   2+u_{3}\ar[ur]&    &      &   &   &    \\
 &&&u_{1}u_{2} \ar[ur]\ar[dr]&&v\ar[ur]&&&&&&&&\\
 &&&&u_{3}\ar[ur]&&&&&&&\\
}$}\\

Using the uni-modular rule for the bottom square in this diagram, we have the equation $vu_{1}u_{2}= h_{0}u_{3}+1$ which becomes $2h_{0}=2$, due to the above expression for $v$, thus $h_{0} = 1$.

\hspace{0.3 cm} Thus we have proved that the quiver $\tilde{Q}$ with variables is equal to $\Lambda'$. This proves that the constructed frieze $\Upsilon$ is associated with $\Lambda'$. $\square$ \\

\begin{rem}
We proved theorem $2.13$ for a seed $\mathcal{G}$ whose associated quiver has two arrows entering the joint (the case of two fork arrows leaving the joint is entirely similar). For a seed $\mathcal{G}_{2}$ whose associated quiver $Q_{2}$ has a fork consisting of one arrow entering and one arrow leaving the joint there are two possibilities.

\hspace{0.3 cm} One possibility is to perform a mutation on one of the fork vertices and to reduce this case to the one considered in the proof (see the discussion following Lemma $2.9$). 

\hspace{0.3 cm} Another possibility is to work directly with the seed $\mathcal{G}_{2} = (Q_{2}, \chi)$. In this case the given proof can be easily modified, namely one has to associate with the vertex labelled by $i=1$ in $Q_{2}'$ the variable $\displaystyle\left( \frac{u_{1}(1+u_{3})}{u_{2}}\right)$ and consider the following underlying graph for $\Lambda'_{2}$ with variables:\\
$\xymatrix @=10pt
{
&
&
&&
&& \\
\text{}
&u_{n} \ar@{-}[r] 
&\ldots \ar@{-}[r] u_{4}
&u_{3} \ar@{-}[r] 
&\displaystyle\left( \frac{u_{1}(1+u_{3})}{u_{2}}\right)  \ar@{-}[r]
&u_{0} \ar@{-}[r]
&\displaystyle\left( \frac{u_{1}(1+u_{3})}{u_{2}}\right)  \ar@{-}[r]
&u_{3} \ar@{-}[r]
&u_{4} \ldots
&u_{n}. \ar@{-}[l]
& 
}
$\\  
\end{rem}

The next step is to establish a formula which allows us to compute directly cluster variables of a cluster algebra of type $\mathbb{D}_{n}$ independently of each other.

\section{Computation of cluster variables: case $\mathbb{D}_{n}$}

Let $\mathcal{G} = (Q, \{u_{1}, ..., u_{n}\})$ be a seed with $Q$ of type $\mathbb{D}_{n}$ with two fork arrows entering or leaving the joint.

\hspace{0.3 cm} In this section, given the initial values $\mathfrak{a}(0, i) = u_{i}, \; i = 1, ...,n$ we compute the values in $\bar{\mathcal{F}}$ by an explicit formula using matrix product (Theorem $3.6$). It is well known, see [ARS10] and also [AD11] of Assem-Dupont, that if $\mathfrak{a}(0, i) = u_{i}$ with $i \in (\mathbb{D}_{n})_{0}$ then all the cluster variables of cluster algebra with the initial seed $\mathcal{G}$ are contained in the part $\mathcal{F}$. Therefore the values of the frieze at the points of the modelled quiver $\bar{\mathcal{F}}$ associated with the seed $\mathcal{G}$ are either cluster variables of the cluster algebra $\mathcal{A}(\mathcal{G})$ or products of two such cluster variables.

\hspace{0.3 cm} According to theorem $2.13$, the modelled quiver $\bar{\mathcal{F}}$ associated with $\mathcal{G}$ is contained in the frieze $\Upsilon$ of type $\mathbb{A}_{2n-1}$ associated with $\Lambda'$. Therefore the computation of values lying in $\bar{\mathcal{F}}$ can be done by considering these values as those in the frieze $\Upsilon$ of type $\mathbb{A}_{2n-1}$ associated with $\Lambda'$. The frieze $\Upsilon$ corresponds to a triangulation without internal triangles, thus we use the result of [ARS10] where the values of a frieze of type $\mathbb{A}_{n}$ were computed.

\begin{defi}
We call boundary a sequence $ c_{1} x_{1} c_{2} x_{2} ... c_{m-1} x_{m-1} c_{m} $  with $x_{i} \in \{ x,y \}$ and $c_{i} \in \mathbb{Q}(u_{1}, u_{2},...,u_{n})$, $i \in \mathbb{Z}$.
\end{defi}

For the variables $\alpha , \beta$ in $\mathbb{Q}(u_{1}, u_{2},..., u_{n})$, replacing the arrows of the form $\alpha\rightarrow \beta$ by $\alpha x \beta $ and those of the form $\alpha \leftarrow \beta$ by $\alpha y \beta$, we can regard $\Lambda'$ as a boundary (by abuse of notation we denote the resulting boundary by $\Lambda'$ as well):\\
$\Lambda'= u_{n}x_{1}u_{n-1}x_{2}...u_{3}x_{n-2}u_{1}u_{2}x_{n-1}1x_{n}u_{1}u_{2}x_{n+1}u_{3}x_{n+2}u_{4}...x_{2n-3}u_{n-1}x_{2n-2}u_{n}$, \\with $x_{k} \in \{x, y \}$.  \\

We give an example showing the boundary $\Lambda'$ for a given seed $\mathcal{G}$.\\
\begin{exem}
 $\xymatrix @=10pt{
&{1}\ar[dr]
&
&
&\text{}\\
\text{For the following quiver $\mathbb{D}_{4}$}\quad
&
&{3}\ar[r]
&{4}
&&\text{}\\  
&{2}\ar[ur]
&
&
}  $
$\xymatrix @=10pt{
&
&
&\text{}\\
& \text{we have} \; \; \textbf{$\bar{\omega}$} :\;{1}\ar[r]
&{3}\ar[r]
&4
&&\text{}\\  
&
&
&
&
}  $

$\xymatrix @=10pt{
&
&
&\text{}\\
& \text{} \; \; \textbf{${Q'}$} :\; {4}
&{3}\ar[l]
&{1}\ar[l]\ar[r]
&{0}\ar[r]
&{1}\ar[r]
&{3}\ar[r]
&{4} \quad \text{and}\quad \Lambda' = u_{4}yu_{3}yu_{1}u_{2}x1xu_{1}u_{2}xu_{3}xu_{4}.
&&\text{}\\  
&
&
&
&
}  $\\
 
\end{exem}

\hspace{0.3 cm} Before giving a formula which allows us to compute each variable lying at the point $(u, v)$ in $\bar{\mathcal{F}}$, we will embed our frieze $\Upsilon$ into the Euclidean plane by rotating it by $45^{\circ}$ clockwise. More precisely, the arrows directed north-east will become horizontal segments and the arrows directed south-east will become vertical segments.

\hspace{0.3 cm} Each boundary may be embedded into the Euclidean plane in the following way: before embedding, a boundary is extended on its end points by $y$ on the left and $x$ on the right. Then $x$ (or $y$) determine the horizontal (or vertical, respectively) segments of a discrete path, that is $x$ (or $y$) corresponds to a segment of the form [$(u , v)\; ,\; (u+1 , v)]$  (or $[(u , v)\; ,\; (u , v+1)]$, respectively) in the plane. The variables $c_{i}$ become thus labels of the vertices of the discrete path. \\

 In Example $3.5$ we give an example of an embedded boundary.\\
 
\hspace{0.3 cm} Let us define the notions of the \textit{ position boundary} $F$ and the word associated with a point in a frieze. \\
 
\begin{defi}
Let $\Lambda'$ be the boundary defined above. We call position boundary $F$ a boundary associated with $\Lambda'$, embedded in the discrete plane. In other words, $F = y \Lambda' x$.\\
\end{defi}

As $\bar{\omega}_{*}$ bounds the modelled quiver $\bar{\mathcal{F}}$ on the left and on the right, the part of the frieze $\Upsilon$ containing $\bar{\mathcal{F}}$ is bordered on the left by $\Lambda'$ and on the right by $^{t}\Lambda'$. Considered in the Euclidean plane, $\bar{\mathcal{F}}$ is contained in the part of the frieze $\Upsilon$ bordered by the position boundary $F = y \Lambda' x$ on one side and $^{t}F = x\; ^{t}\Lambda' y$ on the other side, where $^{t}F$ is the transpose of $F$. Note that there is a correspondence induced by the transposition between the points of the boundaries $F$ and $^{t}F$.

\hspace{0.3 cm} From now on, $\bar{\mathcal{F}}$ and $\Upsilon$ stand for the modelled quiver and the frieze $\Upsilon$ considered in the Euclidean plane, that is after replacing north-east and south-east arrows by horizontal and vertical segments respectively.

\hspace{0.3 cm} Consider a point $(u,v)$ in $\bar{\mathcal{F}}$ and its horizontal and vertical projections on the boundaries $F$ and $^{t}F$. Considering the following three situations suffices to define the word associated with the point $(u , v)$:\\
- The two projections lie on the boundary $F$,\\
- The two projections lie on the boundary $ ^{t}F$,\\
- The horizontal projection lies on $F$ and the vertical projection lies on $^{t}F$.

\hspace{0.3 cm} Recall that in the proof of theorem $2.13$ the frieze $\Upsilon$ was constructed by completing $\bar{\mathcal{F}}$ downward using the uni-modular rule. In the Euclidean plane, after the clockwise rotation by $45^{\circ}$, this implies that $\bar{\mathcal{F}}$ coincides with the upper right part of the frieze $\Upsilon$ lying above the oblique line joining the variables $u_{1}u_{2}$ (see Example $3.5$). Then, the three situations are sufficient because they represent the three regions of $\bar{\mathcal{F}}$ illustrated in Example $3.5$ (recall that our aim is to compute the values lying at the points in $\bar{\mathcal{F}}$).

\hspace{0.3 cm} In each situation, we define the word associated with the point $(u,v)$ which we will use to compute the variable lying at the point $(u,v)$ in $\bar{\mathcal{F}}$.\\
 
 \begin{defi}
 Let $(u,v)$ be a point in $\bar{\mathcal{F}}$ and let $F$ be the position boundary. The word associated with the point $(u,v)$ is the portion of the boundaries $F$ or $^{t}F$ determined as follows:\\
 - if the horizontal and vertical projections of the point $(u,v)$ lie on the boundary $F$ then the word associated with the point $(u,v)$ is a portion of $F$ delimited by these projections.\\
 -  if the horizontal and vertical projections of the point $(u,v)$ lie on the boundary $^{t}F$ then the word associated with the point $(u,v)$ is a portion of $^{t}F$ delimited by these projections with $x$ and $y$ interchanged.\\
 - if the horizontal projection of the point $(u, v)$ lies on the boundary $F $ and its vertical projection lies on $ ^{t}F$, the word associated with the point $(u, v)$ is the portion of the boundary $F$ delimited by the horizontal projection and the vertex of $F$ corresponding to the vertical projection under the transposition.
 
 \end{defi}
 
The following example shows how to associate a word with a point in each situation.\\
 
  \begin{exem}

    Consider the following boundaries:  \\$F =  yu_{6}y u_{5}y u_{4}xu_{3}xu_{1}u_{2}xxu_{1}u_{2}yu_{3}yu_{4}xu_{5}xu_{6}x $ and its transpose \\ $^{t}F =  xu_{6}yu_{5}yu_{4}xu_{3}xu_{1}u_{2}yyu_{1}u_{2}yu_{3}yu_{4}xu_{5}x u_{6}y $ \; corresponding to the quiver 

 $
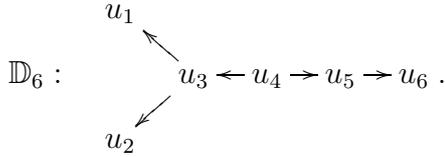
\xymatrix@=10pt{
&u_{1}
&
&
&\text{}\\
 \mathbb{D}_{6} :
&
&u_{3}\ar[ul]\ar[dl]
&u_{4}\ar[r]\ar[l]
&u_{5}\ar[r]
&u_{6} \; .
&&\text{ }\\  
&u_{2}
&
&
}  $\\

  The frieze of type $\mathbb{A}_{2n-1}$ associated with $\Lambda'$, embedded in the discrete plane, is as follows:
   
$
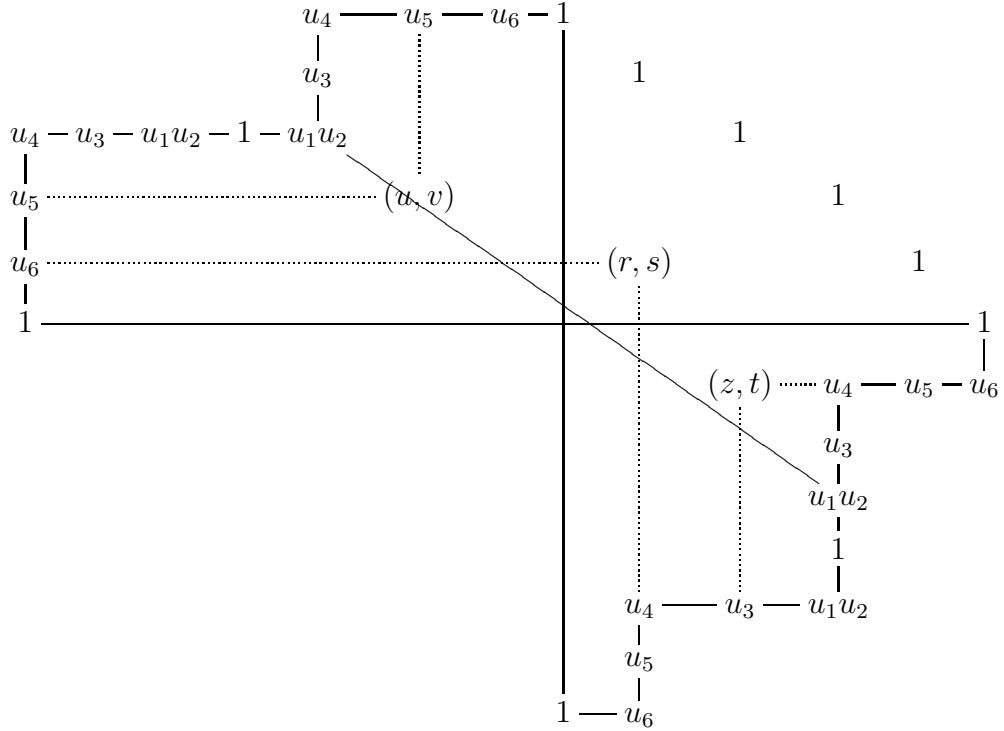
\xymatrix@=7pt{
&      &      &   &    &&  &   u_{4}\ar@{-}[d]\ar@{-}[r]   &     u_{5}\ar@{-}[r]\ar@{.}[ddd]  &   u_{6}\ar@{-}[r]   & 1 \ar@{-}[dddddddddddd] &         \\
&  &   &   &    &&  & u_{3}\ar@{-}[d]     &   &       &    &1  &   \\
&   &  & u_{4}\ar@{-}[r]\ar@{-}[d] &  u_{3}\ar@{-}[r]&u_{1}u_{2}\ar@{-}[r]& 1\ar@{-}[r]&  u_{1}u_{2} \ar@{-}[rrddrrrrdddd] &   &  && &1 &  \\
&    &  & u_{5}\ar@{-}[d]\ar@{.}[rrrrr] &  & &  &&    (u,v)&  & &      &&1&   & \\
&    &  & u_{6}\ar@{-}[d]\ar@{.}[rrrrrrrr]  &   &   &  &&  &&&(r,s)\ar@{.}[dddddd] && &    1  &   & \\
&    &  & 1\ar@{-}[rrrrrrrrrrrr]  &   &   &  &&   &&&&& &&    1\ar@{-}[d]  &   & \\
   && &  &  &   &  & &    &&  &  &(z,t)\ar@{.}[dddd]\ar@{.}[r]&u_{4}\ar@{-}[r]\ar@{-}[d] &u_{5}\ar@{-}[r]&   u_{6} &   & \\
&&    &  &  &   &   &  && & &&& u_{3}\ar@{-}[d] &      &   & \\
&&    &  &  &   &   &  &&&  &&& u_{1}u_{2}\ar@{-}[d] &      &   & \\
&&    &  &  &   &   &  && & &&& 1\ar@{-}[d] &      &   & \\
&& &   &  &  &   &  &&& &  u_{4}\ar@{-}[r]\ar@{-}[d]  &u_{3}\ar@{-}[r]&u_{1}u_{2}  &      &   & \\
&& &&   & && &  & &  &u_{5}\ar@{-}[d]&&  &      &   & \\
&& &&&   &&& & & 1\ar@{-}[r] &  u_{6} &   &    &&&&  &      &   & \\
}$\\

\vspace{0.5 cm}

\hspace{0.3 cm}  We illustrate the three situations from definition $3.4$ by: \\
\begin{enumerate}
\item The projections of the point $M=(u,v)$ lie on the boundary $F$, the word associated with this point $M=(u,v)$ is $u_{5}yu_{4}xu_{3}xu_{1}u_{2}xxu_{1}u_{2}yu_{3}yu_{4}xu_{5} $.\\

\item The projections of the point $L=(z,t)$ lie on the boundary $ ^{t}F$, thus we are looking at the portion $u_{3}xu_{1}u_{2}y1yu_{1}u_{2}yu_{3}yu_{4}$ of the boundary and the word associated with this point $L=(z,t)$ is $u_{3}yu_{1}u_{2}x1xu_{1}u_{2}xu_{3}xu_{4}$.\\

\item The horizontal and vertical projections of the point $N=(r,s)$ lie respectively on the boundaries $F$ and $^{t}F$, the word associated with this point $N=(r,s)$ is a portion of the boundary $F$ delimited by the variable $u_{6}$ corresponding to the horizontal projection and the variable $u_{4}$ corresponding to the vertical projection by transposition. Therefore the word associated with the point $N=(r,s)$ is : $ u_{6}yu_{5}yu_{4}xu_{3}xu_{1}u_{2}xxu_{1}u_{2}yu_{3}yu_{4} $.\\ 
\end{enumerate}
\end{exem}

  For $a , b \in  \mathbb{Q}(u_{1}, u_{2},...,u_{n})$ we define the following matrices:\\

 $ M(a , x , b)= \left( \begin{array}{cccccc}
&a &  &1 &  \\
&0 &  &b& \\
\end{array} \right) ;\quad  M(a , y , b)= \left( \begin{array}{cccccc}
&b &  &0 &  \\
&1 &  &a& \\
\end{array} \right).$ \qquad \qquad  $(2)$ \\

 The next theorem allows us to compute the variables in  $\bar{\mathcal{F}}$ independently of each other.

\begin{theo}
Let $\bar{\mathcal{F}}$ be a modelled quiver associated with a seed  $\mathcal{G}$ of type $\mathbb{D}_{n}$ with arrows of the fork both entering or both leaving the joint. Consider a point $(u,v)$ in $\bar{\mathcal{F}}$ with the following associated word $b_{0}x_{1}b_{1} x_{2}...b_{n}x_{n+1}b_{n+1}$, $n \geq 1$, $x_{i} \in \{ x,y \}$, $ b_{i} \in \mathbb{Q}(u_{1}, u_{2},...,u_{n}) $. The value at this point $(u,v)$ is given by the function $\textbf{T}: \mathbb{Z}^2 \rightarrow \mathbb{Q}(u_{1}, u_{2},...,u_{n})$ defined as follows: \\
\begin{enumerate}
\item[a.] For each point $(u, v)$ in $\bar{\mathcal{F}}$ lying on the boundary $F$ or $^{t}F$, $\textbf{T}(u, v)$ coincides with the value of the frieze $\Upsilon$ at this point.

\item[b.] If $(u,v)$ is a point such that its two projections lie on one of the boundaries $F$ or $^{t}F$, then

$\textbf{T}(u,v)= \displaystyle\frac{1}{b_{1}b_{2} . . . b_{n}}(1,b_{0})\prod_{i=2}^{n}M(b_{i-1},x_{i},b_{i})\left( \begin{array}{cccccc}
1 \\
b_{n+1} \\
       
\end{array} \right),$ \qquad \qquad \qquad \qquad $(3)$

\item[c.] If $(u,v)$ is a point such that its horizontal projection lies on $F$ and its vertical projection lies on $^{t}F$, then 

$\textbf{T}(u,v)= \displaystyle\frac{1}{b_{1}b_{2} . . . b_{n}}(1,b_{0})\prod_{i=2}^{n}M(b_{i-1},x_{i},b_{i})\left( \begin{array}{cccccc}
b_{n+1} \\
1 \\       
\end{array} \right).$ \qquad \qquad \qquad \qquad $(4)$
\\
\end{enumerate}
\end{theo}

The proof of theorem $3.6$ will rely on the identification of the values in $\bar{\mathcal{F}}$ with the values in the frieze $\Upsilon$ of type $\mathbb{A}_{2n-1}$ of theorem $2.13$. We need the following two lemmas from [ARS10], where $\mathbb{K}$ denotes some field. We recall the proofs for the convenience of the reader.\\

\begin{lem}
\begin{enumerate}
\item Let $A \in \mathcal{M}_{2\times 2}(\mathbb{K})$ be a square matrix, $\; \lambda , \lambda' \in \mathcal{M}_{1\times 2}(\mathbb{K}),$ two row matrices, $\;  \gamma ,\gamma' \in \mathcal{M}_{2\times 1}(\mathbb{K})$ two column matrices.\\
Consider the scalars given by the following matrix products:\\ $p = \lambda A \gamma ,\quad q = \lambda A \gamma' , \quad r = \lambda' A \gamma ,\quad s = \lambda' A \gamma'$.\\ Then $  det\left( \begin{array}{cccccccccc}
&p &  &q &  \\
& r &  &s& \\
\end{array} \right)$ = $detA$ \; $  det \left( \begin{array}{cccccc}
\lambda \\
\lambda' \\
       
\end{array} \right) \; det ( \gamma , \gamma' ).$
\item  Let $a, b_{1}, b_{2},...,b_{k}, b \in \mathbb{K},$  and consider the following row matrices $\mathcal{M}_{1\times 2}(\mathbb{K})$: \\ $\lambda' = (1 , a ) M(b_{1}, x ,b_{2}). . . M(b_{k-1}, x , b_{k})M(b_{k},y,b)$ and   $\lambda = (1, b_{k})$ with $M \in \mathcal{M}_{2\times 2}(\mathbb{K})$ given by $(2)$ . Then $det \left( \begin{array}{cccccc}
\lambda' \\
\lambda \\
       
\end{array} \right) = b_{1} b_{2} . . .b_{k}b .$\\

\end{enumerate}
\end{lem}

\textbf{Proof}
\begin{enumerate}
\item We have  $\left( \begin{array}{cccccccccc}
&p &  &q &  \\
& r &  &s& \\
\end{array} \right)=  \left( \begin{array}{cccccccccc}
&\lambda A \gamma &  &\lambda A \gamma' &  \\
& \lambda' A \gamma &  &\lambda' A \gamma'& \\
\end{array} \right) = \left( \begin{array}{cccccc}
\lambda \\
\lambda' \\
       
\end{array} \right) A ( \gamma , \gamma' )$ and the result follows.\\
\item Let $N =  M(b_{1}, x ,b_{2}). . . M(b_{k-1}, x , b_{k})$. A direct computation gives $$\qquad \qquad N = \left( \begin{array}{cccccccccc}
&  b_{1} b_{2} . . .b_{k-1}&  &u &  \\
& 0 &  & b_{2} . . .b_{k}& \\
\end{array} \right) \qquad \qquad (5)$$ with some $u \in \mathbb{K}. $\\

Using the equality $(1,a)N = ( b_{1} b_{2} . . .b_{k} , u +   ab_{2} . . .b_{k})$ we obtain by a straightforward computation the value of $\lambda'$: \\
$\lambda' = (1 , a )  N \left( \begin{array}{cccccccccc}
&b &  &0 &  \\
& 1 &  &b_{k}& \\
\end{array} \right)=\displaystyle \left[ ( b_{1} b_{2} . . .b_{k-1})b+u+ ab_{2} . . .b_{k} , (u+  ab_{2} . . .b_{k})b_{k} \right] .\\$
It follows that \;$ det  \left( \begin{array}{cccccc}
\lambda' \\
\lambda \\
\end{array} \right) 
= - (u+  ab_{2} . . .b_{k})b_{k} + (bb_{1} . . .b_{k-1}+ u+  ab_{2} . . .b_{k})b_{k}. $ \nonumber\\ 
\end{enumerate}

By developing and simplifying we get the result. {$\square$}  \\

\vspace{1cm}

\begin{lem}

Let $A \in \mathcal{M}_{2\times 2}(\mathbb{K}), a , b_{1}, . . . b_{k}, b, c, c_{1}, . . . ,c_{l},d \in \mathbb{K},\quad k,l\geq 1$.\\ Consider the real numbers $p,q,r,s$ defined by the following matrix products:\\
\begin{eqnarray}
p &=& (1, b_{k})A \left( \begin{array}{cccccc}
1 \\
c_{1} \\
       
\end{array} \right), \nonumber \\
q &=&(1 , b_{k} )A M(c, x ,c_{1}) M(c_{1}, y , c_{2}) . . . M(c_{l-1},y,c_{l})\left( \begin{array}{cccccc}
1 \\
d \\
       
\end{array} \right), \nonumber\\
r &=&(1 , a )M(b_{1}, x ,b_{2}) . . . M(b_{k-1}, x , b_{k})M(b_{k},y,b) A \left( \begin{array}{cccccc}
1 \\
c_{1} \\
       
\end{array} \right), \nonumber\\
s &=& (1 , a )M(b_{1}, x ,b_{2}) . . . M(b_{k-1}, x , b_{k})M(b_{k},y,b) A M(c,x,c_{1}) \times \nonumber\\ &\times & M(c_{1},y,c_{2}) . . . M(c_{l-1},y,c_{l})\left( \begin{array}{cccccc}
1 \\
d \\
\end{array} \right) \nonumber
\end{eqnarray}
with $M \in \mathcal{M}_{2\times 2}(\mathbb{K})$ given by $(2)$.
 Then 
   $ det \left( \begin{array}{cccccccccc}
&p &  &q &  \\
& r &  &s& \\
\end{array} \right) = b_{1} . . . b_{k}bcc_{1} . . . c_{l}\; det A$. \\

\end{lem}

\textbf{Proof:} Let $\gamma = \left( \begin{array}{cccccc}
1 \\
c_{1} \\
       
\end{array} \right)  \in \mathcal{M}_{2\times 1}(\mathbb{K}), \; \; \lambda = (1 , b_{k} ), \\ \lambda' = (1 , a )M(b_{1}, x ,b_{2}) . . . M(b_{k-1}, x , b_{k})M(b_{k},y,b) \in \mathcal{M}_{1\times 2}(\mathbb{K})$ and \\
$\gamma' =  M(c, x ,c_{1}) M(c_{1}, y , c_{2}) . . . M(c_{l-1},y,c_{l})\left( \begin{array}{cccccc}
1 \\
d \\
       
\end{array} \right) \in \mathcal{M}_{2\times 1}(\mathbb{K}) $.\\
Applying lemma $3.7$-$1.$ we have $ det \left( \begin{array}{cccccccccc}
&p &  &q &  \\
& r &  &s& \\
\end{array} \right) =  det$ $A$ $ det \left( \begin{array}{cccccc}
\lambda \\
\lambda' \\
       
\end{array} \right) det( \gamma , \gamma' )$, also by lemma $3.7$-$2.$ we have  $ det \left( \begin{array}{cccccc}
\lambda' \\
\lambda \\
       
\end{array} \right)$ =$ - b_{1}. . . b_{k}b$ \; \; and symmetrically \\ $ det(\gamma , \gamma' )$ = $- cc_{1}. . . c_{l}$ . Thus $ det \left( \begin{array}{cccccccccc}
&p &  &q &  \\
& r &  &s& \\
\end{array} \right)$ = $b_{1} . . . b_{k}bcc_{1} . . . c_{l}$ $ det A$. {$\square$}  \\

In the proof of theorem $3.6$, we are going to apply the above two lemmas with $\mathbb{K}=\mathbb{Q}(u_{1}, u_{2},...,u_{n})$.\\

\textbf{Proof of theorem 3.6:} Due to theorem $2.13$, $\bar{\mathcal{F}}$ is contained in the frieze of type $\mathbb{A}_{2n-1}$ associated with $\Lambda'$. Values in this frieze are obtained, starting from the initial values of $\Lambda'$, by applying the uni-modular rule (see proof of theorem $2.13$). Note that by definition, the function $T$ takes the initial values of $\Lambda'$ on the boundary. Therefore to prove the theorem it suffices to show that the function $T$ satisfies the uni-modular rule: for four adjacent points forming a square in the discrete plane $\mathbb{Z}^{2}$, the images satisfy the uni-modular rule, that is\\
 $$ det \left( \begin{array}{ccccc}
T(u,v+1) &  &T(u+1,v+1)   \\
 T(u,v) &  &T(u +1,v) \\
\end{array} \right) = 1.$$\\

Let us prove the uni-modular rule for the case when all four points have both projections lying on $F$. In this case formula $(3)$ is applied to calculate the four corresponding values. The case when either formula $(4)$ or both formulas $(3)$ and $(4)$ are applied to calculate the four values are entirely similar due to the fact that the points $(u, v+1)$ and $(u, v)$ or $(u+1, v+1)$ and $(u+1, v)$ have the same vertical projection.\\
       
       Consider the following scheme which represents four adjacent points $(u,v) ,\; \\ (u, v+1), \; (u+1, v+1), \; (u+1, v)$ forming a square in the frieze.\\

$
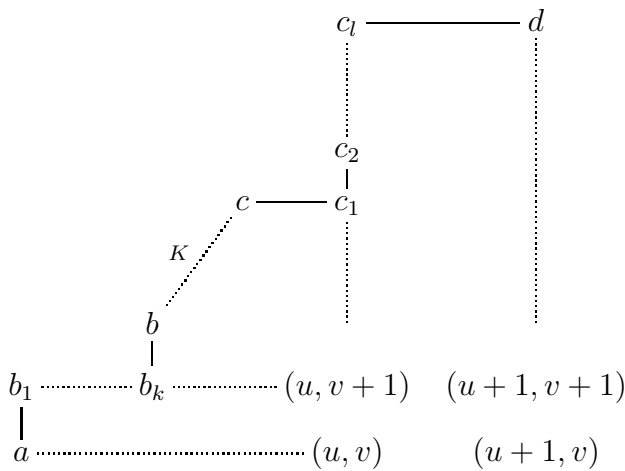
\xymatrix@=7pt{
&&& &   &   &   &  c_{l}\ar@{.}[ddd]&d\ar@{.}[ddddddd]\ar@{-}[l]&      &               \\
 &   &   &  &  &    &     &   &   &       \\
 &&  &  &  &&   &  &  &   &          \\
&&  & & && &   c_{2}& &          \\
&&& & & &   c\ar@{-}[r] & c_{1}\ar@{.}[ddd]\ar@{-}[u]&   &          \\
& &  & &  &  &   &   &          \\
& &  & &  &  &   &   &          \\
& &  & & b\ar@{.}[uuurr]^{K} &  &   &   &          \\
&b_{1}\ar@{.}[rrr]& & &b_{k}\ar@{-}[u]\ar@{.}[rrr]&&&(u,v+1)&(u+1,v+1)&&\\
 &a\ar@{.}[rrrrrr]\ar@{-}[u]& &   &  &    &   & (u,v)&(u+1,v)&    &     \\
}$\\

Let us find the words associated with the four given points.\\

Denote by $K = x_{1}\alpha_{1}x_{2} . . . . .\alpha_{n-1}x_{n}$ the portion of the boundary between $b$ and $c$ as illustrated in the figure, $\alpha_{i} \in \mathbb{K}$ and $x_{i} \in \{x, y \}$.\\
Then the words associated with the points $(u,v) ,\; (u, v+1), \; (u+1, v+1)$ \\ and $(u+1, v)$ \; are respectively $ayb_{1}x. . . .b_{k-1}xb_{k}ybKcxc_{1}; \; b_{k}ybKcxc_{1}; \\ b_{k}ybKcxc_{1}yc_{2} . . .  .c_{l-1}yc_{l}xd$ and $ayb_{1}x. . . .b_{k-1}xb_{k}ybKcxc_{1}yc_{2} . . .  .c_{l-1}yc_{l}xd$.
Let a matrix $A$ be constructed from $K$ as follows: $$A = M(b, x_{1}, \alpha_{1})M(\alpha_{1}, x_{2}, \alpha_{2}) . . . . M(\alpha_{n-1}, x_{n}, c) \qquad \qquad (6)$$  with $M \in \mathcal{M}_{2 \times 2}(\mathbb{K})$ given by $(2)$ \; and let a scalar $D$ be obtained from $K$ as follows: \; $$D = b\alpha_{1}. . . . \alpha_{n-1}c \in \mathbb{K} .\qquad \qquad \qquad \qquad \qquad \qquad \qquad \qquad (7)$$\\
We obtain the following equalities by using the definition of the map $\textbf{T}$:\\
 $T(u,v)=\displaystyle\frac{1}{b_{1}. . . b_{k} D}(1 , a )M(b_{1}, x ,b_{2}) . . . M(b_{k-1}, x , b_{k})M(b_{k},y,b) A \left( \begin{array}{cccccc}
1 \\
c_{1}\end{array} \right),$ \\
$T(u+1,v)=\displaystyle\frac{1}{b_{1}. .b_{k} D c_{1}. .c_{l}} (1 , a )M(b_{1}, x ,b_{2}) . .M(b_{k-1}, x , b_{k})M(b_{k},y,b)\times \\  A M(c, x ,c_{1}) M(c_{1}, y , c_{2}) . .M(c_{l-1},y,c_{l})\left( \begin{array}{cccccc}
1 \\
d \\
\end{array} \right)$,\\
$T(u+1,v+1)=\displaystyle\frac{1}{ D c_{1}. . . c_{l}}(1 , b_{k} )A M(c, x ,c_{1}) M(c_{1}, y , c_{2}) . . . M(c_{l-1},y,c_{l})\left( \begin{array}{cccccc}
1 \\
d \\
       
\end{array} \right)$,\\
$T(u,v+1) =\displaystyle \frac{1}{D} (1, b_{k})A \left( \begin{array}{cccccc}
1 \\
c_{1} \\
       
\end{array} \right)$.\\
Let us define $p , q , r , s $ as in lemma $3.8$. Then the above formulas can be rewritten as follows:\\
$T(u,v)=\displaystyle\frac{r}{b_{1}. . . b_{k} D}$, \quad $T(u+1,v)=\displaystyle\frac{s}{b_{1}. . . b_{k} D c_{1}. . . c_{l}}$, \quad \\
 $T(u+1,v+1)=\displaystyle\frac{q}{ D c_{1}. . . c_{l}}$,\quad $T(u,v+1) =\displaystyle \frac{p}{D}$.\\

Then we have
 \begin{eqnarray}
T(u,v+1)T(u+1, v) - T(u,v)T(u+1,v+1) &=& \displaystyle\frac{1}{b_{1}. . . b_{k} D^{2} c_{1}. . . c_{l}} (ps-qr) \nonumber \\
&=& \displaystyle\frac{b_{1} . . .b_{k}bcc_{1} . . . c_{l}}{b_{1}. . . b_{k} D^{2} c_{1}. . . c_{l}}\; \rm det A . \nonumber 
\end{eqnarray}

From definition $(6)$ of the matrix $A$ we find its determinant: \\$detA = b\alpha_{1}^{2}\alpha_{2}^{2} . . . \alpha_{n-1}^{2}c$. Using expression $(7)$ for $D$ we have the following result:\\
$ T(u,v+1)T(u+1, v) - T(u,v)T(u+1,v+1) = 1$.\\

\hspace{0.3 cm} It remains to prove that $T$ satisfies the uni-modular rule at the boundary.\\

\hspace{0.3 cm} Suppose that one of the four points lies on the boundary, for example $(u, v+1) = b_{k}$. We have the following scheme:\\

$
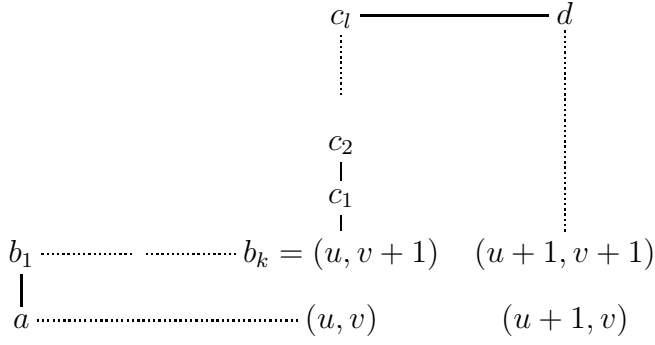
\xymatrix@=7pt{
&&& &   &   &   &  c_{l}\ar@{.}[dd]&d\ar@{.}[ddddd]\ar@{-}[l]&      &               \\
 &   &   &  &  &    &     &   &   &       \\
 &&  &  &  &&   &  &  &   &          \\
&&  & & && &   c_{2}& &          \\
&&& & & &    & c_{1}\ar@{-}[d]\ar@{-}[u]&   &          \\
&b_{1}\ar@{.}[rrr]& & &\ar@{.}[rrr]&&&b_{k} = (u,v+1)&(u+1,v+1)&&\\
 &a\ar@{.}[rrrrrr]\ar@{-}[u]& &   &  &    &   & (u,v)&(u+1,v)&    &     \\
}$\\

\vspace{0.5 cm}

The words associated with the points $(u,v) ,\;  (u+1, v+1)$ and $(u+1, v)$ \\ are respectively $ayb_{1}x. . . .b_{k-1}xb_{k}; \;  b_{k}yc_{1}yc_{2} . . .  .c_{l-1}yc_{l}xd$ and \\ $ayb_{1}x. . . .b_{k-1}xb_{k}yc_{1}yc_{2} . . .  .c_{l-1}yc_{l}xd$.

\hspace{0.3 cm} We obtain the following equalities by using the definition of the map $\textbf{T}$:\\
 $T(u,v)=\displaystyle\frac{1}{b_{1}. . . b_{k-1}}(1 , a )M(b_{1}, x ,b_{2}) . . . M(b_{k-2}, x , b_{k-1}) \left( \begin{array}{cccccc}
1 \\
b_{k}\end{array} \right),$ \\
$T(u+1,v)=\displaystyle\frac{1}{b_{1}. .b_{k} c_{1}. .c_{l}} (1 , a )M(b_{1}, x ,b_{2}) . .M(b_{k-1}, x , b_{k})M(b_{k},y,c_{1})\times \\ M(c_{1}, y , c_{2}) . .M(c_{l-1},y,c_{l})\left( \begin{array}{cccccc}
1 \\
d \\
\end{array} \right)$,\\
$T(u+1,v+1)=\displaystyle\frac{1}{ c_{1}. . . c_{l}}(1 , b_{k} ) M(c_{1}, y , c_{2}) . . . M(c_{l-1},y,c_{l})\left( \begin{array}{cccccc}
1 \\
d \\
       
\end{array} \right)$,\\
$T(u,v+1) = b_{k}$.\\

\hspace{0.3 cm} Let us first compute the products of matrices appearing in $T(u, v)$ and $T(u+1, v+1)$.

\hspace{0.3 cm}Let $W =  M(b_{1}, x ,b_{2}). . . M(b_{k-2}, x , b_{k-1})$ be the matrix product from $T(u, v)$ corresponding to the boundary between $b_{1}$ and $b_{k}$. Since all segments in this portion are horizontal, we compute $W$ as in $(5)$: 
 $$W = \left( \begin{array}{cccccccccc}
  b_{1} b_{2} . . .b_{k-2}&  &v  \\
 0 &  & b_{2} . . .b_{k-1} \\
\end{array} \right), \; \text{with some}\; v \in \mathbb{K}.$$\\

Let now $Z =  M(c_{1}, y ,c_{2}). . . M(c_{l-1}, y , c_{l})$ be the matrix product from $T(u+1, v+1)$ corresponding to the boundary between $c_{1}$ and $c_{l}$. Since all segments in this portion are vertical, we compute $Z$ analogously to $(5)$:
$$ Z = \left( \begin{array}{cccccccccc}
&   c_{2} . . .c_{l}&  &0 &  \\
& z &  & c_{1} . . .c_{l-1}& \\
\end{array} \right), \; \text{with some}\; z \in \mathbb{K}.$$\\ 

\hspace{0.3 cm}Using the equalities $ (1,a)W = ( b_{1} b_{2} . . .b_{k-2} , v +   ab_{2} . . .b_{k-1})$ and\\

 $ Z\left( \begin{array}{cccccc}
1 \\
d \\
\end{array} \right) = \left( \begin{array}{cccccc}
c_{2}...c_{l} \\
z + dc_{1}...c_{l-1} \\
\end{array} \right)   $ we obtain the following values by putting $(\tilde{a}, \tilde{b}) = (1,a)W$ and $\left( \begin{array}{cccccc}
\tilde{c} \\
\tilde{d} \\
\end{array} \right) = Z\left( \begin{array}{cccccc}
1 \\
d \\
\end{array} \right)$:\\
 $T(u,v)=\displaystyle\frac{1}{b_{1}. . . b_{k-1}}(\tilde{a} , \tilde{b} ) \left( \begin{array}{cccccc}
1 \\
b_{k}\end{array} \right),$ \\
$T(u+1,v)=\displaystyle\frac{1}{b_{1}. .b_{k} c_{1}. .c_{l}} (\tilde{a} , \tilde{b} )\left( \begin{array}{cccccccccc}
  1+  c_{1}b_{k-l}&  &b_{k}   \\
 b_{k} &  & b_{k}^{2} \\
\end{array} \right)\left( \begin{array}{cccccc}
\tilde{c} \\
\tilde{d} \\
\end{array} \right)$,\\
$T(u+1,v+1)=\displaystyle\frac{1}{ c_{1}. . . c_{l}} (1, b_{k}) \left( \begin{array}{cccccc}
\tilde{c} \\
\tilde{d} \\
       
\end{array} \right)$,\\
$T(u,v+1) = b_{k}$.\\

Then we have by a straightforward computation\\

 \begin{eqnarray}
T(u,v+1)T(u+1, v) - T(u,v)T(u+1,v+1) &=& \displaystyle\frac{1}{b_{1}. . . b_{k-1} c_{1}. . . c_{l}} (\tilde{a}\tilde{c}c_{1}b_{k-1}) \nonumber \\
&=& 1 . \nonumber 
\end{eqnarray}

We obtain the same result with other situations where two or three of the four points lie on the boundary. $\square$\\

\begin{rem}

The values in $\bar{\mathcal{F}}$ are the cluster variables of cluster algebra of type $\mathbb{D}_{n}$ except those lying on the diagonal line passing through the variable $u_{1}u_{2}$. On this line, the values are products of two cluster variables (products created by the passage from ${\mathcal{F}}$ to $\bar{\mathcal{F}}$). Note that the pairs of cluster variables forming these products (in the case of both arrows entering or leaving the joint) are given by fractions whose numerators are equal and denominators coincide up to exchanging $u_{1}$ and $u_{2}$, which appear in denominators with exponent one (see [BMR09-1]). Therefore, a product of these two variables is a perfect square divided by $u_{1}u_{2}$. Thus each value $\textbf{T}(u,v) = U \times V$ on the diagonal in question gives rise to two cluster variables $U$ and $V$ as follows:\\

 $\textbf{T}(u,v)= U \times V$ where \\
 $U = \displaystyle\frac{1}{u_{1}}\left[ \frac{u_{1}u_{2}}{b_{1}b_{2} . . . b_{n}}(1,b_{0})\prod_{i=2}^{n}M(b_{i-1},x_{i},b_{i})\left( \begin{array}{cccccc}
1 \\
b_{n+1} \\
       
\end{array} \right)\right]^{\dfrac{1}{2}}$ and \\
$V = \displaystyle\frac{1}{u_{2}}\left[ \frac{u_{1}u_{2}}{b_{1}b_{2} . . . b_{n}}(1,b_{0})\prod_{i=2}^{n}M(b_{i-1},x_{i},b_{i})\left( \begin{array}{cccccc}
1 \\
b_{n+1} \\
       
\end{array} \right)\right]^{\dfrac{1}{2}}$. 
\end{rem}

We give now an example to illustrate the above results.\\

\begin{exem}

 $\xymatrix @=10pt{
&{1}\ar[dr]
&\text{}\\
\text{For the following quiver $\mathbb{D}_{4}$: }
&
&{3}\ar[r]
&{4}
&&\text{}\\  
&{2}\ar[ur]
&
&
}  $
$\xymatrix @=10pt{
&
&
&\text{}\\
& \text{we have} \; \; \textbf{$\bar{\omega}$} :\;{1}\ar[r]
&{3}\ar[r]
&4
&&\text{}\\  
&
&
&
&
}  $

$\xymatrix @=10pt{
&
&
&\text{}\\
& \text{} \; \; \textbf{${Q'}$} :\; {4}
&{3}\ar[l]
&{1}\ar[l]\ar[r]
&{0}\ar[r]
&{1}\ar[r]
&{3}\ar[r]
&{4} \quad \text{and}\quad \Lambda' = u_{4}yu_{3}yu_{1}u_{2}x1xu_{1}u_{2}xu_{3}xu_{4}.
&&\text{}\\  
&
&
&
&
}  $\\

The boundary associated with $\Lambda'$ is $F  = yu_{4}yu_{3}yu_{1}u_{2}x1xu_{1}u_{2}xu_{3}xu_{4}x$ and the modelled quiver $\bar{\mathcal{F}}$ bordered by $F$ and $^{t}F$ embedded in the plane gives the following diagram:\\
  
$
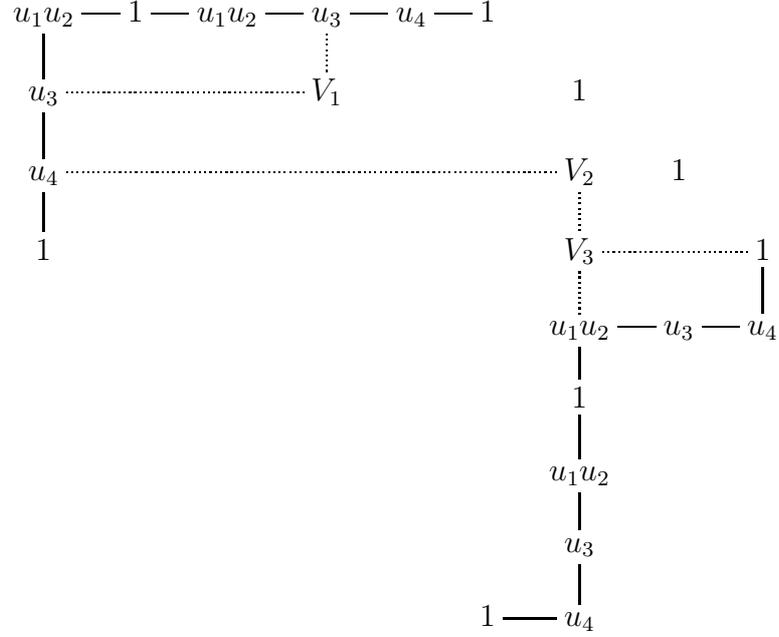
\xymatrix@=14pt{
&      &      &   &    & u_{1}u_{2}\ar@{-}[r]\ar@{-}[d]& 1\ar@{-}[r] & u_{1}u_{2}\ar@{-}[r]   &     u_{3}\ar@{-}[r]\ar@{.}[d]  &   u_{4}\ar@{-}[r]   & 1 &         \\
&  &   &   &    &u_{3}\ar@{-}[d]&  &     & V_{1} \ar@{.}[lll] &       &    &1  &   \\
&   &  & & &u_{4}\ar@{-}[d]& &   &   &  &&V_{2} \ar@{.}[d]\ar@{.}[llllll] &1 &  \\
&    &  &  &  & 1&  && &  & &  V_{3}\ar@{.}[d] \ar@{.}[rr]    &&1 \ar@{-}[d]&   & \\
&    &  &  &   &   &  &&  &&& u_{1}u_{2}\ar@{-}[d]\ar@{-}[r]&u_{3}& u_{4}\ar@{-}[l] &      &   & \\
&    &  &  &   &   &  &&   &&&1\ar@{-}[d]&& &&   &   & \\
   && &  &  &   &  & &     &&  &  u_{1}u_{2}\ar@{-}[d] &&&&  &   & \\
&&    &  &  &   &   &  && & & u_{3}\ar@{-}[d] &&&      &   & \\
&&    &  &  &   &   &  &&& 1\ar@{-}[r]  &u_{4}&&  &      &   & \\
&&    &  &  &   &   &  && & &&&  &      &   & \\
}$\\

Where $V_{1}$ represents a product of two cluster variables on the bottom line in $\bar{\mathcal{F}}$, $V_{2}$ and $V_{3}$ represent some cluster variables on a line in $\bar{\mathcal{F}}$ except the bottom line.\\

We are going to compute the cluster variables in these three positions.\\

The words associated with $V_{1}$, $V_{2}$ and $V_{3}$ are respectively $u_{3}yu_{1}u_{2}x1xu_{1}u_{2}xu_{3}$, \\
$u_{4}yu_{3}yu_{1}u_{2}$ and $u_{1}u_{2}yu_{3}yu_{4}x1.$\\

Then we have, by applying theorem $3.6$, the following results:\\
\begin{eqnarray}
V_{1} &=& \displaystyle\frac{1}{u_{1}^{2}u_{2}^{2}} \left( 1,u_{3}\right)  \left( \begin{array}{cccccc}
&u_{1}u_{2} &  &1 &  \\
&0 &  &1& \\
\end{array} \right) \left( \begin{array}{cccccc}
&1 &  &1 &  \\
&0 &  &u_{1}u_{2}& \\
\end{array} \right) \left( \begin{array}{cccccc}
1 \\
u_{3} \\     
\end{array} \right) \nonumber \\
&=&\displaystyle\frac{(1+u_{3})^{2}}{u_{1}u_{2}}  \nonumber\\
&=&\displaystyle\frac{1+u_{3}}{u_{1}}\times \frac{1+u_{3}}{u_{2}}.\nonumber
\end{eqnarray}
$V_{1}$ being placed on the bottom line in $\bar{\mathcal{F}}$, we write this value in the form of a product of two cluster variables as in remark $3.9$. The other two positions correspond to the following cluster variables of cluster algebra of type $\mathbb{D}_{4}$: 
\begin{eqnarray}
V_{2}&=& \displaystyle\frac{1}{u_{3}} \left( 1,u_{4}\right)  \left( \begin{array}{cccccc}
u_{1}u_{2} \\
1 \\     
\end{array} \right) \nonumber\\
&=& \displaystyle\frac{u_{4}+u_{1}u_{2}}{u_{3}}, \nonumber
\end{eqnarray}
and
\begin{eqnarray}
V_{3}&=& \displaystyle\frac{1}{u_{3}u_{4}} \left( 1,u_{1}u_{2}\right)  \left( \begin{array}{cccccc}
&u_{4} &  &0 &  \\
&1 &  &u_{3}& \\
\end{array} \right) \left( \begin{array}{cccccc}
1 \\
1 \\     
\end{array} \right) \nonumber \\
&=&\displaystyle\frac{u_{4}+u_{1}u_{2}(1+u_{3})}{u_{3}u_{4}}.  \nonumber
\end{eqnarray}\\
\end{exem}

ACKNOWLEDGEMENTS. The author is grateful to Ibrahim Assem and Vasilisa Shramchenko for useful discussions and careful reading of the manuscript and also to the department of mathematics of the University of Sherbrooke where this work was done. The author gratefully acknowledges support from CIDA (Canadian International Development Agency) fellowship.

{Kodjo Essonana MAGNANI\\
   D\'epartement de math\'ematiques\\
   Universit\'e de Sherbrooke\\
2500, boul. de l'Universit\'e,\\
Sherbrooke, Qu\'ebec,
J1K 2R1\\
   Canada\\
   Kodjo.essonana.magnani@USherbrooke.ca}

\end{document}